\documentclass[11pt,a4paper,twoside]{article}
\usepackage{bm, amsmath, amssymb, amsthm} 

\def\calf{{\cal F}}
\def\<{\langle}
\def\>{\rangle}
\def\eps{\varepsilon}

\newcommand\const{\operatorname{const}}
\newcommand\tr{\operatorname{Tr}}
\newcommand\Div{\operatorname{div}}
\newcommand\id{\operatorname{id}}

\def\vol{\operatorname{vol}}

\def\eq{\hspace*{-1.5mm}&=&\hspace*{-1.5mm}}
\def\plus{\hspace*{-1.5mm}&+&\hspace*{-1.5mm}}
\def\minus{\hspace*{-1.5mm}&-&\hspace*{-1.5mm}}

\newtheorem{example}{Example}

\newtheorem{lem}{Lemma}
\newtheorem{prop}{Proposition}
\newtheorem{thm}{Theorem}

\author{Vladimir Rovenski\footnote{Mathematical Department, University of Haifa, Mount Carmel, 31905 Haifa,  Israel
        \newline e-mail: {\tt vrovenski@univ.haifa.ac.il} }
 }

\title{The Reeb type formula \\ for codimension-one foliated $(\alpha,\beta)$-spaces}

\begin{document}

\date{}

\maketitle


\begin{abstract}
Integral formulae for foliated Riemannian manifolds provide obstructions for existence of foliations or compact leaves of them
with given geometric properties.
Recently, we associated a new Riemannian metric to a codimension-one foliated Finsler space
 and proved integral formulae for general and for Randers spaces.
In the paper, we study this metric for a wider class of codimension-one foliated $(\alpha,\beta)$-spaces
and embody it in a set of metrics that can be viewed as a perturbed metric associated with $\alpha$.
For such metrics we calculate the Weingarten operator of the leaves and deduce the Reeb type integral formula,
which can be used for $(\alpha,\beta)$-spaces, e.g. Randers and Kropina spaces.


\end{abstract}

\section*{Introduction}

Recent decades brought increasing interest in Finsler spaces, especially, in extrinsic geometry of their hypersurfaces,
see \cite{sh2}.
The theory of $(\alpha,\beta)$-metrics has been developed into a fruitful branch of Finsler geometry,
see~e.g.~\cite{shsh2016}--\cite{yz}.
A~{Finsler structure} $F$ on a manifold $M$ is a family of Minkowski norms in tangent spaces $T_pM$
depending smoothly on a point $p\in M$.
The~$(\alpha,\beta)$-\textit{norms}, introduced in 1972 by M.\,Matsumoto as extension of Randers and Kropina norms, 
form a special class of `computable' Minkowski norms, see e.g. \cite{ma92}.

Being the central topic of the extrinsic geometry of foliations, {\it integral formulae}
(i.e., integral relations for invariants of the shape operator of the leaves and Riemannian curvature)
provide obstructions for existence of foliations with given geometric properties
and have applications in different areas of geometry and analysis, see \cite{lw2,rw2}.
The~first known integral formula by G.\,Reeb~\cite{re} tells~us~that the~total mean curvature $H$ of the leaves
of a codimension-1 foliated closed Riemannian manifold is zero; thus, either $H\equiv0$ or $H(p)H(q)<0$  for some points $p\ne q$.
This was extended in \cite{rw2} into infinite series of integral formulas with the higher order mean curvatures,
see survey in \cite{r25,rw1}.

 In \cite{rw3,rw4}, to a codimension-one foliated Finsler space $(M,F)$
we associated a new Riemannian metric $g$ and derived its Riemann curvature in terms of~$F$.
For Randers metric $F=\alpha+\beta$ we derived the shape operator of the leaves and
obtained new integral formulae, comparing the Reeb integral formula for $g$ and the metric associated with $\alpha$.
 Continuing this line of research, the paper presents new integral formulae for a
 foliated Finsler space with $(\alpha,\beta)$-metric and such examples as Randers and Kropina metrics.

 Our main results (verified with assistance of Maple) is Theorem~\ref{T-1-1}
 generalizing Reeb integral formula for a manifold with a general $(\alpha,\beta)$-metric.
 In Section~\ref{sec:ab}, we survey necessary facts and prove auxiliary lemmas for $g$ considered as a perturbed source metric.
In Section~\ref{sec:shape}, we calculate the shape operator of the leaves for $g$,
then find its trace in Section~\ref{sec:Reeb} to prove the Reeb type integral formula.
 Results of the paper are valid for foliations and 1-forms defined outside a finite union
 of closed submanifolds of codimension $\ge 2$ under convergence of some integrals, see discussion in \cite{lw2}.
The singular case is important since there exist a number of manifolds which admit no (smooth) codimension-one foliations,
while all of them admit such foliations and non-singular 1-forms $\beta$ outside some ``set of singularities".

\section{The $(\alpha,\beta)$-metric}
\label{sec:ab}

A \textit{Minkowski norm} on a vector space $V^{m+1}$ is a function $F:V^{m+1}\to[0,\infty)$ with the  properties of regularity, positive 1-homogeneity and strong convexity, see \cite{shsh2016}:

\smallskip
M$_1:$ $F\in C^\infty(V^{m+1}\setminus \{0\})$,

M$_2:$ $F(\lambda\,y)=\lambda F(y)$ for all $\lambda>0$ and $y\in V^{m+1}$,

M$_3:$ For any $y\in V^{m+1}\setminus \{0\}$, the following symmetric bilinear form is positive definite:
\begin{equation}\label{E-acsiom-M3}
 g_y(u,v)=\frac12\,\frac{\partial^2}{\partial s\,\partial t}\,\big[F^2(y+su+tv)\big]_{|\,s=t=0}\,.
\end{equation}

\noindent
By (M$_2$)\,--\,(M$_3$), $g_{\lambda y}=g_{y}\ (\lambda>0)$ and $g_y(y,y)=F^2(y)$.

 Let $\phi:(-b_0,b_0)\to(0,\infty)$ be a smooth function.
Suppose that $a(\cdot\,,\cdot)=\<\cdot\,,\cdot\>$ is a scalar product on $V^{m+1}$,
$\alpha(y)=\sqrt{\<y,y\>}$ the Euclidean norm,
and $\beta$ a linear form on $V^{m+1}$ of the norm $b=\alpha(\beta)<b_0$.
 The $(\alpha,\beta)$-\textit{metric} is defined on $V^{m+1}\setminus\{0\}$ by \cite{ma92,shsh2016}
\begin{equation}\label{E-ab-def}
  F(y)=\alpha(y)\,\phi(s),\quad s=\beta(y)/\alpha(y).
\end{equation}

For better understanding the concept of $(\alpha,\beta)$-metric, recall the following.

\begin{lem}[\cite{shsh2016}]
 The function \eqref{E-ab-def} is a Minkowski norm on $V^{m+1}$ for any $\alpha$ and $\beta$ with
 $\alpha(\beta) <b_0$ if and only if the function $\phi(s)$ satisfies
\begin{equation}\label{E-phi-cond}
 \phi-s\,\phi\,' +(b^2-s^2)\,\phi\,'' > 0,
\end{equation}
where $s$ and $b$ are arbitrary numbers with $|s|\le b<b_0$.

Taking $b=s$ in \eqref{E-phi-cond}, one can see that $\phi$ obeys $\phi-s\,\phi\,'>0$ when $|s|<b_0$.
\end{lem}

For $(\alpha,\beta)$-metric, the bilinear form $g_y\ (\alpha(y)=1)$ in \eqref{E-acsiom-M3}
(which can be viewed as a perturbed scalar product $\<\cdot\,,\cdot\>$)
and the ratio $\sigma_g(y)=\det g_y/\det a$ are given by
\begin{eqnarray}\label{E-c-value0}
 g_y(u,v)\eq \rho\<u,v\> +\rho_0\beta(u)\beta(v) +\rho_1(\beta(u)\<y,v\> +\beta(v)\<y,u\>) -\rho_1\beta(y)\<y,u\>\<y,v\>,\\
\label{E-F001e}
\nonumber
 \sigma_g(y) \eq \rho^{m-1}\big(
 \rho^2 +\rho_0\rho_1\beta(y)^3+\rho_1^2\beta(y)^2+(\rho-\rho_0b^2)\rho_1\beta(y)+(\rho\rho_0-\rho_1^2)b^2\,\big)\\
 \eq \phi^{m+2}(\phi-s\,\phi\,')^{m-1}[\,\phi-s\,\phi\,'+(b^2-s^2)\phi\,''\,] >0,
\end{eqnarray}
see e.g. \cite{shsh2016}, where $\rho,\rho_0,\rho_1$ are the following functions of $s$:
\begin{equation*}
 \rho = \phi(\phi-s\,\phi\,'),\quad
 \rho_0 = \phi\,\phi\,''+(\phi\,')^2,\quad
 \rho_1 = \rho' = -s(\phi\,\phi\,''+(\phi\,')^2)+\phi\phi\,'.
\end{equation*}

 Let $W\subset V^{m+1}$ be a hyperplane and $N$ a unit normal to $W$ with respect to $\<\cdot\,,\cdot\>$, i.e.,
\[
 \<N,v\>=0\quad (v\in W),\qquad {\<N,N\>} = 1.
\]
The 'musical isomorphisms' $\sharp$ and $\flat$ will be used for rank one and symmetric rank 2 tensors.
For~example, $\<\beta^\sharp,u\>=\beta(u)=u^\flat(\beta^\sharp)$ for $u\in V^{m+1}$.
Let $\beta^{\sharp\top}$ be the projection of $\beta^\sharp$ onto $W$.

For any Minkowski norm, there are exactly two normal directions to $W$, see \cite{sh2}, which are opposite when $F$ is {reversible}.
Hence, there is a unique $\alpha$-unit vector $n\in V^{m+1}$, which is $g_n$-orthogonal to $W$
and belongs to the same half-space as $N$, i.e.,
\[
 g_n(n,v)=0\quad (v\in W),\quad \<n,n\>=1,\quad \<n,N\> > 0.
\]
Then
$\nu=F(n)^{-1}n$ is an $F$-unit normal to~$W$, where $F(n)=\phi(s)$ and $s=\beta(n)$. 
Set $g:=g_n$, see \eqref{E-c-value0} with $y=n$, and assume
\begin{equation}\label{E-gamma2-discr}
 (b^2-\beta(N)^2)(\rho_1+\rho_0\beta(n))^2\le\rho^2 \ \Longleftrightarrow\
 (b^2-\beta(N)^2)^{1/2}\,|\phi'(s)|\le \phi(s)-s\,\phi'(s).
\end{equation}
Note that $b^2-\beta(N)^2\ge0$. By \eqref{E-F001e}, $\det g = \sigma_g(n)\det a$ for $m>1$.
Define the quantities
\begin{eqnarray*}
 && \gamma_1 :=({\rho_1+\rho_0\beta(n)})/{\rho} ={\phi\,'}({\phi-s\,\phi\,'})_{\,|s=\beta(n)},\\
 && \gamma_2:= \rho_0-\gamma_1\rho_1(\beta(n)\gamma_1+2)
 ={\phi\,(\phi^2\phi\,''-\phi(\phi\,')^2+s(\phi\,')^3)}/{(\phi-s\,\phi\,')^2}_{\,|s=\beta(n)},
\end{eqnarray*}
and set
 $\hat c=\gamma_1\beta(N) +\sqrt{1 -\gamma_1^2(b^2-\beta(N)^2)}$
-- by \eqref{E-gamma2-discr}, the discriminant is nonnegative.

\begin{lem}\label{L-c-value}
Let \eqref{E-gamma2-discr} holds, then
\begin{eqnarray}\label{E-c-value}
 && n = \hat c\,N-\gamma_1\beta^\sharp,\\
\label{E-c-value2}
 && g(u,v) =  \rho\,\<u,v\> +\gamma_2\beta(u)\,\beta(v)\quad (u,v\in W)\,, \\
\label{E-c-value2c}
 && g(n,n) = \rho + \rho_1\beta(n)+\rho_0\beta(n)^2 = \phi(\beta(n))^2.
\end{eqnarray}
\end{lem}

\proof
From \eqref{E-c-value0} with $y=n$, $u=n$ and $v\in W$ and $g(n,v)=0$ we find
\begin{equation}\label{E-gnv}
 \rho\,\< n +\gamma_1 \beta^\sharp,\ v\> =0\quad(v\in W).
\end{equation}
From \eqref{E-gnv} and $\rho>0$ we conclude that $n+\gamma_1\beta^{\sharp}=\hat c\,N$ for some real $\hat c$.
Using
\[
 1 = \<n,n\> = \hat c^{\,2}-2\,\hat c\,\beta(N)\gamma_1 +(b\,\gamma_1)^2
\]
and $b^2=\alpha(\beta^{\sharp\bot})^2+\alpha(\beta^{\sharp\top})^2$, we get two real solutions
\[
 \hat c_{1,2} =\gamma_1\beta(N)\pm\sqrt{1 -\gamma_1^2(b^2-\beta(N)^2)}\,,
\]
the~greater value (with $+$) provides the required $\<n,N\> > 0$, that proves \eqref{E-c-value}.
Thus,
\[
 \beta(n)=\beta(\hat c\,N-\gamma_1\beta^\sharp)=\hat c\,\beta(N) - \gamma_1b^2.
\]
Finally, \eqref{E-c-value2} follows from \eqref{E-c-value0}, \eqref{E-c-value} and $\<n,u\>=-\gamma_1\beta(u)\ (u\in W)$,
and \eqref{E-c-value2c} follows from \eqref{E-c-value0}.
\qed

\smallskip

Assume that $\gamma_2$ is sufficiently `small' (relative to $\rho>0$), i.e.,
\begin{equation}\label{E-gamma3-cond}
 \rho+(b^2-\beta(N)^2)\,\gamma_2 \ne 0,
\end{equation}
and define the quantity
\[
 \gamma_3=-\frac{\gamma_2}{\rho+(b^2-\beta(N)^2)\,\gamma_2}\,.
\]

The following lemma relates metrics $g$ and $\<\cdot\,,\cdot\>$ along $W$.

\begin{lem}\label{L-zZ}
Let \eqref{E-gamma3-cond} holds, $u,U\in W$ and
\[
 g(u,v)=\<U,v\>\quad (\forall\,v\in W),
\]
then
\begin{equation}\label{E-u-U}
 \rho\,u = U+\gamma_3\beta(U)\,\beta^{\sharp\top} .
\end{equation}
\end{lem}

\proof
By \eqref{E-c-value2} we have
\[
 g(u,v) = \<\rho\,u +\gamma_2\beta(u)\beta^{\sharp},\,v\>.
\]
By conditions, since $u,U$ and $\beta^{\sharp\top}$ belong to~$W$, we find
 $\rho\,u +\gamma_2\beta(u)\beta^{\sharp\top} = U$.
Applying $\beta$, yields
\[
 \beta(u)(\rho+(b^2-\beta(N)^2)\gamma_2)=\beta(U)
\]
and then \eqref{E-u-U}.
\qed

\begin{example}\label{Ex-rand1}\rm
For $\phi(s)=1$ we have Riemannian metric $F=\alpha$.
Then for a hyperplane $W\subset V^{m+1}$ and the metric $g=g_n=\<\cdot\,,\cdot\>$
we obtain
 $n=N$, $\hat c=1=\rho$, and $\rho_0=\rho_1=\gamma_1=\gamma_2=\gamma_3=0$.

Recently, some progress was achieved about many particular cases of $(\alpha,\beta)$-metrics.

(i)~For $\phi(s)=1+s,\ |s|\le b<b_0=1$, we have the metric $F=\alpha+\beta$,
first introduced by a physicist G.\,Randers
to consider the unified field theory.
For a hyperplane $W\subset V^{m+1}$ and the metric $g=g_n$ with $s=\beta(n)$ we obtain, see also \cite{rw4},
\[
 n=\hat c\,N-\beta^\sharp,\quad
 \beta(n)=c\,\hat c-1,\quad
 \|n\|_g = c\,\hat c,
\]
where $\hat c = c+\beta(N)$ and $c=\sqrt{1-(b^2-\beta(N)^2)}\in(0,1]$. Then
\begin{eqnarray*}
 && \rho=c\,\hat c,\quad \rho_0=\rho_1=\gamma_1=1,\quad
 \gamma_2=-c\,\hat c,\quad
 \gamma_3=c^{-2}.
\end{eqnarray*}
Conditions \eqref{E-gamma2-discr} and \eqref{E-gamma3-cond} become trivial (i.e., $c>0$).
Next, $\sigma_g(n)=(1+\beta(n))^{\,m+2}$ and
\begin{equation*}
 g(u,v) = (1+\beta(n))\<u,v\> - \beta(n)\<n,u\>\,\<n,v\> +\beta(u)\,\<n,v\> +\beta(v)\,\<n,u\> +\beta(u)\,\beta(v).
\end{equation*}

(ii)~A subclass of $(\alpha,\beta)$-metrics having the form $F = \alpha^{l+1}/\beta^l\ (l\ne0)$, i.e., $\phi(s) = 1/s^l\ (s>0)$,
are called \textit{generalized Kropina metrics}, see \cite{ma92}.
The \textit{Kropina metric}, i.e., $l=1$,
first introduced by L.\,Berwald, was investigated by V.K.\,Kropina in 1961.
Although there are singularities, the metric has applications in general dynamical systems.
Then \eqref{E-phi-cond} reads $s^2+b^2>0$ and we have
\[
 \rho=2/s^2,\quad \rho_0=3/s^4,\quad \rho_1=-4/s^3.
\]
For a hyperplane $W\subset V^{m+1}$ and the metric $g=g_n$ with $s=\beta(n)$ we obtain
\begin{eqnarray*}
 && \gamma_1=-1/(2\beta(n)),\quad \gamma_2=\gamma_3=0,\\
 && \beta(n) = {\rm sign}(\beta(N))\sqrt{b(b+|\beta(N)|)/2},\\
 && \hat c=\sqrt{1+(b^2-\beta(N)^2)/(4\beta(n)^2)}-\beta(N)/(2\beta(n)).
\end{eqnarray*}
Hence, $g(u,v) = \rho\,\<u,v\>\ (u,v\in W)$.
Condition \eqref{E-gamma3-cond} becomes trivial, and \eqref{E-gamma2-discr} becomes
$4\,\beta(n)^2 \ge b^2-\beta(N)^2$ and is also satisfied.
\end{example}

In next sections we deal with general $(\alpha,\beta)$-metrics under conditions \eqref{E-gamma2-discr}, \eqref{E-gamma3-cond}.

\section{The shape operator of $\calf$ and the curvature of normal flow}
\label{sec:shape}

Let $F$ be a Finsler metric on a connected smooth manifold $M^{m+1}\ (m\ge2)$ such that $F(p, y)$ is
a \textit{general $(\alpha,\beta)$-metric} on the tangent space $T_pM$ for any $p\in M$, see~\cite{yz},
i.e., $\phi=\phi(p,s)$, $\alpha=\alpha(p,y)$ and $\beta=\beta(p,y)$ in \eqref{E-ab-def} depend also on a point $p\in M$.
Let $\bar\nabla$ be the Levi-Civita connection of Riemannian metric $\<\cdot\,,\cdot\>$ associated with $\alpha$.

Given a transversally oriented codimension-one foliation $\calf$ of $(M,F)$,
there exists a glo\-bally defined $F$-normal to the leaves smooth vector field $n$,
which defines a Riemannian metric $g:=g_n$ (see \eqref{E-c-value0} for $y=n$), with the Levi-Civita connection $\nabla$.
Due to Section~\ref{sec:ab}, let $N$ be a unit $a$-normal vector field to $\calf$ and the following hold:
\[
 \<n,N\> >0,\quad \<n,n\>=1.
\]
Then $s=\beta(n)$ is a function on $M$ and $\nu=n/F(n)$ is a $g$-unit normal (and an $F$-unit normal). 
The shape operators $\bar A,A^g:T\calf\to T\calf$ of the leaves are defined by
\begin{equation}\label{E-Ag}
 \bar A(u)=-{\bar\nabla_u\,N},\quad
 A^g(u)=-{\nabla_u\,\nu}\quad (u\in T\calf).
\end{equation}
The curvature vectors of $\nu$- and $N$- flows for metrics $g$ and $\langle\cdot,\cdot\rangle$ on $M$ are, respectively,
\begin{equation*}
 Z=\nabla_\nu\,\nu,\quad \bar Z=\bar\nabla_N\,N.
\end{equation*}
The $(1,1)$-tensors $\bar\nabla u:TM\to TM$ and its $\langle\cdot,\cdot\rangle$-conjugate $(\bar\nabla u)^t: TM\to TM$
are defined~by
\[
 (\bar\nabla u)\,(v)=\bar\nabla_v\,u,\quad
 \<(\bar\nabla\,u)^t(v),w\>=\<v,(\bar\nabla\,u)(w)\>\quad (v,w\in TM).
\]
The {deformation tensor},
 $\overline{\rm Def}_u =\frac12\,(\bar\nabla u+(\bar\nabla u)^t)$,
measures the degree to which the flow of a vector field $u$ distorts
$\langle\cdot,\cdot\rangle$.
We use the~same notation $\overline{\rm Def}_u$ for its $\langle\cdot,\cdot\rangle$-conjugate $(1,1)$-tensor.
For a vector field $X$ and 1-form $\omega$ on $M$, set ${\rm Sym}(X\otimes\omega)=\frac12\,(X\otimes\omega + X^\flat\otimes\omega^\sharp)$.

Remark that a parallel vector field $\beta^\sharp$ ($\bar\nabla\beta^\sharp=0$)
forms a constant angle with the leaves of $\calf$ if and only if $\beta(N)=\const$ and $b=\const$.

\begin{prop}\label{L-Dx}
Let $\calf$ be a codimension-one foliation of $M^{m+1}\ (m\ge2)$ with a general $(\alpha,\beta)$-metric $F$ and conditions \eqref{E-gamma2-discr}, \eqref{E-gamma3-cond}. Then
\[
 \rho\,\|n\|_g A^g =-{\cal A}-\gamma_3\,(\beta\circ{\cal A})\otimes\beta^{\sharp\top},
\]
where the linear operator ${\cal A}: T\calf\to T\calf$ is given by
\begin{equation}\label{E-Dx2}
 {\cal A} = -\rho\,\hat c\bar A -\rho\,\gamma_1(\overline{\rm Def}_{\beta^\sharp})^\top
 + \frac12\,n(\rho)\id^\top +\,{\rm Sym}(U\otimes\beta^\top),
\end{equation}
and
\begin{eqnarray}\label{E-U}
\nonumber
 U \eq \frac12\,n(\gamma_2)\,\beta^{\sharp\top} +\gamma_2\bar\nabla^\top_n\,\beta^{\sharp\top} -\rho\,\bar\nabla^\top\gamma_1
  +\hat c\,\rho_1(1+\beta(n)\gamma_1)\,( (\hat c-\beta(N)\gamma_1)\bar Z +\gamma_1\bar A(\beta^{\sharp\top}))\\
 \plus (\rho_0-\rho_1\gamma_1)\,\big(\beta(N)\bar\nabla^\top\hat c
 -\,(\gamma_1/2)\bar\nabla^\top b^2 -b^2\bar\nabla^\top\gamma_1 +(\hat c -\beta(N)\gamma_1)(\beta(N)\bar Z -\bar A(\beta^{\sharp\top}))\big).\qquad
\end{eqnarray}
Moreover, if $\,b$ and $\beta(N)$ are constant then
\[
 {\cal A} = -\rho\,\hat c\bar A -\rho\,\gamma_1(\overline{\rm Def}_{\beta^\sharp})^\top +{\rm Sym}(U\otimes\beta^\top),
\]
where
\begin{eqnarray*}
 U \eq \gamma_2\bar\nabla^\top_n\,\beta^{\sharp\top}
 +\hat c\,\rho_1(1+\beta(n)\gamma_1)\,\big((\hat c-\beta(N)\gamma_1)\bar Z +\gamma_1\bar A(\beta^{\sharp\top})\big)\\
 \plus (\rho_0-\rho_1\gamma_1)\,(\hat c -\beta(N)\gamma_1)(\beta(N)\bar Z -\bar A(\beta^{\sharp\top})).
\end{eqnarray*}
\end{prop}

\proof By the formula for the Levi-Civita connection of $g$,
\begin{eqnarray}\label{eqlevicivita}
\nonumber
 2\,g(\nabla_u v, w) \eq u(g(v,w)) + v(g(u,w)) - w(g(u,v)) \\
 \plus g([u, v], w) - g([u, w], v) - g([v, w], u),
\end{eqnarray}
where $u, v, w\in C^\infty(TM)$, and equalities
\[
 g(u,n)=0=g(v,n),\quad g([u,v],n)=0\quad (u,v\in T\calf),
\]
we~have
\begin{equation}\label{E-LC-g}
 2\,g(\nabla_u\,n, v) = n(g(u,v)) +g([u,n],v) +g([v,n],u)\quad (u,v\in T\calf).
\end{equation}
Assume $\bar\nabla_X^\top\,u=\bar\nabla_X^\top\,v=0$ for $X\in T_pM$ at a given point $p\in M$.
Using \eqref{E-c-value0} and \eqref{E-c-value2}, we get
\begin{eqnarray*}
 n(g(u,v)) \eq n\big(\rho\<u,v\> +\gamma_2\beta(u)\,\beta(v)\big) = n(\rho)\<u,v\> +n(\gamma_2)\beta(u)\beta(v)\\
 \plus \gamma_2\big(\beta(u)(\bar\nabla_n(\beta^\top))(v) +\beta(v)(\bar\nabla_n(\beta^\top))(u)\big),\\
 g([u,n],v)\eq\rho\<\bar\nabla_u\,n,v\>+\rho_0\beta([u,n])\,\beta(v)+\rho_1[\beta([u,n])\<n,v\>+\beta(v)\<n,[u,n]\>]\\
 \minus \rho_1\beta(n)\<n,[u,n]\>\<n,v\>
 =-\rho\,\hat c\<\bar A(u),v\> -\rho\gamma_1\<\bar\nabla_u\,\beta^\sharp,v\> -\rho\,\beta(v)\<\bar\nabla\gamma_1,u\>\\
 \plus(\rho_0{-}\rho_1\gamma_1)\big\<\beta(N)\bar\nabla\hat c {-}(\gamma_1/2)\bar\nabla b^2 {-}b^2\bar\nabla\gamma_1
 {+} (\hat c {-}\beta(N)\gamma_1)(\beta(N)\bar Z {-} \bar A(\beta^{\sharp\top})), u\big\>\beta(v)\\
 \plus \hat c\,\rho_1(1+\beta(n)\gamma_1)\<(\hat c-\beta(N)\gamma_1)\bar Z+\gamma_1\bar A(\beta^{\sharp\top}),\,u\>\beta(v).
\end{eqnarray*}
where $u,v\in T\calf$.
Formula for $g([v,n],u)$ is obtained from $g([u,n],v)$ after change $u\leftrightarrow  v$.
Substituting the above into \eqref{E-LC-g}, we find
\[
 g(\nabla_u\,n, v)=\<{\cal A}(u),\,v\>,
\]
where ${\cal A}$ is given in \eqref{E-Dx2}--\eqref{E-U}.
In particular,
\begin{eqnarray*}
 \<2\,{\cal A}(u),\beta^{\sharp\top}\> \eq -2\,\rho\,\hat c\<\bar A(\beta^{\sharp\top}),u\>
  -2\,\rho\gamma_1\<\overline{\rm Def}_{\beta^\sharp}(\beta^{\sharp\top}),u\> \\
 \plus n(\rho)\beta(u) + \beta(u)\beta(U) +U(u)(b^2-\beta(N)^2).
\end{eqnarray*}
Using Lemma~\ref{L-zZ} and
$g(\nabla_u\,n,v) =-\|n\|_g\,g(A^g(u),v)$, see \eqref{E-Ag}, we get
\begin{equation*}
 -2\,\rho\,\|n\|_g A^g(u) =2\,{\cal A}(u)+\gamma_3\,\<2{\cal A}(u),\,\beta^{\sharp\top}\>\,\beta^{\sharp\top}
\end{equation*}
that proves the first claim.
Note that $\beta(n),\rho,\rho_i,\gamma_i,\hat c$ are constant when $b$ and $\beta(N)$ are constant.
From the above the second claim follows.
\qed

\begin{prop}\label{P-ZZ}
Let $\calf$ be a codimension-one foliation of $M^{m+1}\ (m\ge2)$ with a general $(\alpha,\beta)$-metric $F$
and conditions \eqref{E-gamma2-discr}, \eqref{E-gamma3-cond}.
Then
\[
 \rho Z = {\cal Z} +\gamma_3\beta({\cal Z})\,\beta^{\sharp\top},
\]
where
\begin{equation*}
  {\cal Z} = \big[p_1\bar\nabla^\top({{\gamma_1}/{\|n\|_g}}) + p_2\bar\nabla^\top({\hat c}/{\|n\|_g})\big]\|n\|_g^{-1}
 +\big[p_3\bar Z + p_4\bar A(\beta^{\sharp\top}) + p_5 \bar\nabla^\top(\beta(N)) \big]\|n\|_g^{-2}
\end{equation*}
and
\begin{eqnarray*}
 && p_1 = \hat c\left( (4\rho_1\gamma_1-\rho_0+3\rho_1\beta(n)\gamma_1^2){b}^2-\rho +{\hat c}^{\,2}\rho_1\beta(n)\right)\beta(N)
  -\rho_1(2\beta(n)\gamma_1+1)\,{\hat c}^{\,2}{\beta(N)}^2 \\
 && \quad -\,\rho_1(\beta(n)\gamma_1+1){b}^2{\hat c}^{\,2} +\gamma_1(\rho_0-2\gamma_1\rho_1-\gamma_1^2\rho_1\beta(n)){b}^{4} +\gamma_1\rho\,b^2,\\
 && p_2 = (\rho_0-2\rho_1\beta(n)\gamma_1^2 -3\rho_1\gamma_1)\,\hat c\,{\beta(N)}^2 +\big(
 \gamma_1(2\gamma_1\rho_1+\gamma_1^2\rho_1\beta(n)-\rho_0){b}^2 \\
 && \quad +\,\rho_1(2+3\beta(n)\gamma_1)\,{\hat c}^{\,2}
 -\gamma_1\rho\big)\beta(N) -{\hat c}^{\,3}\rho_1\beta(n) +\left(\rho-\gamma_1\rho_1(\beta(n)\gamma_1+1){b}^2\right)\hat c,\\
 && p_3 = \gamma_1(3\gamma_1\rho_1 +2\gamma_1^2\rho_1\beta(n) -\rho_0)\,\hat c\,{\beta(N)}^3
    +\big( (\rho_0-5\rho_1\beta(n)\gamma_1^2 -5\rho_1\gamma_1)\,{\hat c}^{\,2} +\gamma_1^2\rho\\
 && \quad +\,\gamma_1^2(\rho_0-2\gamma_1\rho_1-\gamma_1^2\rho_1\beta(n)){b}^2 \big) {\beta(N)}^2
 +\big( 2\rho_1(1+2\beta(n)\gamma_1)\,{\hat c}^{\,3} +\gamma_1\hat c\,\big((3\gamma_1\rho_1\\
 && \quad +\,2\gamma_1^2\rho_1\beta(n)-\rho_0){b}^2
 -2\rho\big)\big)\beta(N) -{\hat c}^{\,4}\rho_1\beta(n) +(\rho-\gamma_1\rho_1(\beta(n)\gamma_1 +1){b}^2){\hat c}^{\,2},\\
 && p_4 = \gamma_1(\rho_0{-}2\gamma_1^2\rho_1\beta(n)-3\gamma_1\rho_1 )\,\hat c\,{\beta(N)}^2
 +\gamma_1\hat c\,\big((\rho_0 {-}2\gamma_1\rho_1 {-}\gamma_1^2\rho_1\beta(n)) {b}^2 {+}\rho\big)
 +[( 4\rho_1\gamma_1\\
 && \quad -\,\rho_0+3\rho_1\beta(n)\gamma_1^2)\,{\hat c}^{\,2}
 {+}\gamma_1^2(2\gamma_1\rho_1{+}\gamma_1^{2}\rho_1\beta(n)-\rho_0){b}^2{-}\gamma_1^2\rho]\beta(N)
 {-}\rho_1(\beta(n)\gamma_1{+}1)\,{\hat c}^{\,3} ,\\
 && p_5 = {\hat c}^{\,3}\rho_1\beta(n)\gamma_1 -\gamma_1\rho_1(2\beta(n)\gamma_1 +1){\hat c}^{\,2}\beta(N)
 +\gamma_1\hat c\,(\gamma_1\rho_1(1+\gamma_1\beta(n)){b}^2 -\rho).
\end{eqnarray*}
Moreover, if $\beta^\sharp$ is tangent to $\calf$ and $\,b=\const$ then
\begin{eqnarray*}
 {\cal Z} \eq
 \big\{\,{\hat c}^{\,2}\big[\rho -{\hat c}^{\,2}\rho_1\beta(n) -\gamma_1\rho_1(\beta(n)\gamma_1+1)\,{b}^2\big]\bar Z \\
 \plus {\hat c}\,\big[ \gamma_1\rho -\rho_1(\beta(n)\gamma_1+1)\,{\hat c}^{\,2}
 +\gamma_1(\rho_0 -2\,\gamma_1\rho_1 -\gamma_1^2\rho_1\beta(n)) {b}^2 \big]\bar A(\beta^{\sharp}) \big\}\|n\|_g^{-2}.
\end{eqnarray*}
\end{prop}

\proof
Extend $X\in T_p\calf$ at a point $p\in M$ onto a neighborhood of $p$ with the property $(\bar\nabla_Y\,X)^\top=0$ for any $Y\in T_pM$.
By the formula for the Levi-Civita connection, \eqref{eqlevicivita}, we obtain at $p$:
\[
 g(Z,X) = g([X,\nu],\nu).
\]
Then, using
$\nu=\|n\|_g^{-1}(\hat c\,N-\gamma_1\beta^{\sharp})$ and $[X, fY]=X(f)Y+f[X, Y]$ we get
\begin{eqnarray}\label{E-g-terms}
\nonumber
 && g([X,\nu],\nu) = (\hat c/\|n\|_g)\,X(\hat c/\|n\|_g)\,g(N,N)
 -X(\hat c\,\gamma_1/\|n\|_g^2)\,g(N,\beta^\sharp)\\
\nonumber
 && +\,(\gamma_1/\|n\|_g) X(\gamma_1/\|n\|_g)\,g(\beta^\sharp,\beta^\sharp)
 +(\hat c/\|n\|_g)^{2} g([X,N],N) \\
 && -\,(\gamma_1\hat c/\|n\|_g^{2})\,[\,g([X,\beta^\sharp],N) + g([X,N],\beta^\sharp)\,]
 +(\gamma_1/\|n\|_g)^{2} g([X,\beta^\sharp],\beta^\sharp).
\end{eqnarray}
To compute first three terms in \eqref{E-g-terms}, by \eqref{E-c-value0} and Lemma~\ref{L-c-value}, we obtain the equalities
\begin{eqnarray*}
 g(\beta^{\sharp},\beta^{\sharp}) \eq \rho\,b^2+\rho_0 b^4 +2\rho_1 b^2\beta(n) -\rho_1\beta(n)^3,\\
 g(N,\beta^{\sharp}) \eq (\rho +\rho_0 b^2+\rho_1\beta(n))\,\beta(N) +\rho_1(b^2 -\beta(n)^2)\<n,N\>,\\
 g(N,N) \eq \rho +\rho_0\beta(N)^2 +2\,\rho_1\beta(N)\<n,N\> -\rho_1\beta(n)\<n,N\>^2.
\end{eqnarray*}
To compute last four terms in \eqref{E-g-terms}, we will use
\begin{eqnarray*}
 && [X,\beta^\sharp] = [X,\beta^{\sharp\top}] +X(\beta(N))N  +\beta(N)\big(\<Z,X\>\,N -\bar A(X)\big),\\
 && [X,N] = \bar\nabla_X N -\bar\nabla_N X =-\bar A(X) -\<\bar\nabla_N X,\,N\>\,N = \<\bar Z,\,X\>\,N -\bar A(X),
\end{eqnarray*}
and by \eqref{E-c-value0} and Lemma~\ref{L-c-value}, obtain the equalities
\begin{eqnarray*}
 g([X,N],\,\beta^{\sharp}) \eq (\rho+\rho_0 b^2 +\rho_1\beta(n))\<[X,N],\beta^{\sharp}\> +\rho_1(b^2 -\beta(n)^2)\<[X,N],n\>,\\
 g([X,\beta^\sharp],\,\beta^{\sharp}) \eq (\rho+\rho_0 b^2 +\rho_1\beta(n))\<[X,\beta^\sharp],\beta^{\sharp}\>
 +\rho_1(b^2 -\beta(n)^2)\<[X,\beta^\sharp],n\>,\\
 g([X,N],\,N) \eq \rho\<[X,N],\,N\> +(\rho_0\beta(N) +\rho_1\<n,N\>)\<[X,N],\beta^{\sharp}\> \\
 \plus \rho_1(\beta(N) -\beta(n)\<n,N\>)\<[X,N],n\>,\\
 g([X,\beta^\sharp],\,N) \eq \rho\<[X,\beta^\sharp],\,N\> +(\rho_0 \beta(N) +\rho_1\<n,N\>)\<[X,\beta^\sharp],\beta^{\sharp}\> \\
 \plus \rho_1(\beta(N) -\beta(n)\<n,N\>)\<[X,\beta^\sharp],n\>.
\end{eqnarray*}
Thus,
\begin{eqnarray*}
 &&\hskip-6mm g([X,\nu],\nu) = (\hat c/\|n\|_g)\,X(\hat c/\|n\|_g)\,[\rho +\rho_0\beta(N)^2 +2\,\rho_1\beta(N)\<n,N\> -\rho_1\beta(n)\<n,N\>^2]\\
 \minus X(\gamma_1\hat c/\|n\|_g^{2})\,[(\rho +\rho_0 b^2+\rho_1 \beta(n))\,\beta(N) +\rho_1(b^2 -\beta(n)^2)\<n,N\>] \\
 \plus (\gamma_1/\|n\|_g) X(\gamma_1/\|n\|_g)\,[\rho\,b^2+\rho_0 b^4 +2\rho_1 b^2 \beta(n) -\rho_1\beta(n)^3]\\
 \plus (\hat c\|n\|_g)^{2}[\rho\<[X,N],\,N\> +(\rho_0\beta(N) {+}\rho_1\<n,N\>)\beta([X,N])
 +\rho_1(\beta(N) {-}\beta(n)\<n,N\>)\<n,[X,N]\>] \\
 \minus (\gamma_1\hat c/\|n\|_g^{2})[\rho\<[X,\beta^\sharp], N\> {+}(\rho_0 \beta(N) {+}\rho_1\<n,N\>)\beta([X,\beta^\sharp])
 {+}\rho_1(\beta(N) {-}\beta(n)\<n,N\>)\<n,[X,\beta^\sharp]\>] \\
 \plus (\gamma_1\hat c/\|n\|_g^{2})\,[(\rho+\rho_0 b^2 +\rho_1\beta(n))\beta([X,N])
 +\rho_1(b^2 -\beta(n)^2)\<n,[X,N]\>] \\
 \plus (\gamma_1^2/\|n\|_g^{2})\,[(\rho+\rho_0 b^2 +\rho_1\beta(n))\beta([X,\beta^\sharp])
 +\rho_1(b^2 -\beta(n)^2)\<n,[X,\beta^\sharp]\>].
\end{eqnarray*}
Note that
 $\<n,N\>=\hat c-\gamma_1\beta(N)$,
 $\beta(n)=\hat c\,\beta(N)-\gamma_1b^2$,
see \eqref{E-c-value}, and
\begin{eqnarray*}
 \<[X,N],\,N\>\eq \<\bar Z,X\>,\\
 \<[X,N],\beta^{\sharp}\> \eq \<\beta(N)\bar Z -\bar A(\beta^{\sharp\top}),\,X\>,\\
 \<[X,N],\,n\> \eq \hat c\,\<[X,N],\,N\>-\gamma_1\<[X,N],\,\beta^\sharp\>
 =\<(\hat c -\gamma_1\beta(N))\bar Z +\gamma_1\bar A(\beta^{\sharp\top}),\,X\>,\\
 \<[X,\beta^\sharp],\,N\> \eq \<\bar\nabla(\beta(N))+\beta(N)\bar Z,\,X\>, \\
 \<[X,\beta^\sharp],\beta^\sharp\> \eq b\,X(b) -\<\bar\nabla_{\beta^\sharp}X,\beta^\sharp\>
 =\<b\,\bar\nabla b +\beta(N)^2\bar Z -\beta(N)\bar A(\beta^{\sharp\top}) ,X\>,\\
 \<[X,\beta^\sharp],\,n\> \eq \hat c\,\<[X,\beta^\sharp],\,N\> -\gamma_1\<[X,\beta^\sharp],\,\beta^\sharp\> \\
 \eq \<(\hat c\,\beta(N)-\gamma_1\beta(N)^2)\bar Z -\gamma_1b\,\bar\nabla b +\gamma_1\beta(N)\bar A(\beta^{\sharp\top}) ,X\>.
\end{eqnarray*}
Hence, $g(Z,\,X) = \<{\cal Z},\,X\>$.
With the help of Lemma~\ref{L-zZ} and Maple, we complete the proof.
\qed

\begin{example}\rm We continue Example~\ref{Ex-rand1}.

(i)~For Randers metric, the shape operator of $\calf$ w.r.t. $g$
obeys, see Proposition~\ref{L-Dx} and \cite{rw4},
\begin{eqnarray*}
 && c\,A^g = \bar A -\frac12\,\hat c^{\,-2}c^{-1}(\hat c\,N-\beta^{\sharp})(c\,\hat c)\,I_m
 +\hat c^{\,-1}\,(\overline{\rm Def}_{\beta^\sharp})_{\,|T\calf}^\top
 +\,\frac12\,\big( U -\bar A(\beta^{\sharp\top})\,\big)\otimes\beta^\top \\
 &&\hskip-2mm
 +\,\frac12\,c^{-2}\big( \bar A(\beta^{\sharp\top}) -\<\bar  A(\beta^{\sharp\top}),\,\beta^{\sharp\top}\>\,\beta^{\sharp\top}
 +\,2\,\hat c^{\,-1}(\overline{\rm Def}_{\beta^\sharp}\,\beta^{\sharp\top})^\top
 +U +\beta(U)\,\beta^{\sharp\top}\big) \,\!^\flat\otimes\beta^{\sharp\top},
\end{eqnarray*}
where
\[
 U=\hat c^{\,-1}(\bar\nabla_n\,\beta^{\sharp\top})^\top -c\bar Z\quad
 c=\sqrt{1-\alpha(\beta^{\sharp\top})^2}\in(0,1].
\]
 Next, we have $\gamma_3=c^{-2}$ and
${\cal Z}=\bar Z -\hat c^{\,-1}\,\bar\nabla^\top\hat c$, see Proposition~\ref{P-ZZ}.

(ii)~For Kropina metric, if $\,b=\const$ and $\beta(N)=0$ then by Proposition~\ref{P-ZZ},
\begin{eqnarray*}
 Z \eq \frac{\hat c^{\,2}( b^2 +2(1 +2\,\hat c^2)\beta(n)^2 -2 b^2 \beta(n)^6)}{\beta(n)^4 g(n,n)}\,\bar Z \\
 \plus \frac{\hat c\,(4\,\hat c^{\,2} \beta(n)^8 +2 b^2 \beta(n)^6 -(1 + 2\,\hat c^{\,2})\beta(n)^2 -2 b^2)}
 {\beta(n)^5 g(n,n)}\,\bar A(\beta^{\sharp\top}).
\end{eqnarray*}
\end{example}

\section{The Reeb type integral formula}
\label{sec:Reeb}

The Reeb integral formula for (an arbitrary) Riemannian metric $g$ on a closed codimension-one foliated manifold $M$ reads
\begin{eqnarray}\label{eq61-g}
 \int_M (\tr A^g)\,{\rm d}\vol_{\,g} = 0.
\end{eqnarray}
For a general $(\alpha,\beta)$-metric on $M$, the volume forms of $\langle\cdot,\cdot\rangle$ and
$g$ with $\sigma_g$ given in \eqref{E-F001e} obey
\begin{equation}\label{E-F001vol}
 {\rm d}\vol_g = \sigma_g(n)\,{\rm d}\vol_a.
\end{equation}
The next theorem generalizes \eqref{eq61-g} for foliated Finsler manifolds with general $(\alpha,\beta)$-metrics.
In particular case of Randers metric, it was obtained in \cite{rw4}, and for $\beta=0$ it reduces to \eqref{eq61-g} for $a$.

\begin{thm}\label{T-1-1}
Let $\calf$ be a codimension-one foliation of $M^{m+1}$ with a general $(\alpha,\beta)$-metric $F$
and conditions \eqref{E-gamma2-discr}, \eqref{E-gamma3-cond}.
Then the following integral formula holds:
\begin{eqnarray}\label{E-IF1-ab-gen0fin}
\nonumber
 &&\hskip0mm\int_M \Big\{\sigma_g(n)(\rho\,\|n\|_g)^{-1} \Big(\,\rho\,\hat c\,\underline{\tr\bar A} +\rho\gamma_1(\,\beta(\bar Z)-N(\beta(N)))\\
\nonumber
 &&-\,(m-1)\,n(\rho)/2 +(1+(b^2-\beta(N)^2)\gamma_3)
 \big[ \rho\,\beta^{\sharp\top}(\gamma_1) -n(\rho+(b^2-\beta(N)^2)\,\gamma_2)/2 \\
\nonumber
 &&-\,(\rho_0-\rho_1\gamma_1)\big\<\beta(N)\bar\nabla\hat c +(\hat c-\beta(N)\gamma_1)(\beta(N)\bar Z-\bar A(\beta^{\sharp\top}))\\
 && \nonumber
 -(\gamma_1/2)\bar\nabla\,b^2 -b^2\bar\nabla\gamma_1,\,\beta^{\sharp\top}\big\>
 -\hat c\,\rho_1(1+\beta(n)\gamma_1)\<(\hat c-\beta(N)\gamma_1)\bar Z +\gamma_1\bar A(\beta^{\sharp\top}),\,\beta^{\sharp}\>\big]\\
 && \nonumber -\,\gamma_2(1+(b^2-\beta(N)^2)\gamma_3)
 \big[ (\gamma_1/2)\beta^\sharp(b^2-\beta(N)^2) +(\hat c-\beta(N)\,\gamma_1)\<\bar A(\beta^{\sharp\top}),\,\beta^{\sharp}\> \big]\Big) \\
 && -\beta^\sharp(\sigma_g(n)\,\gamma_1/\|n\|_g)\Big\}\,{\rm d}\vol_a = 0.
\end{eqnarray}
Moreover, if $\,b$ and $\beta(N)$ are constant then
\eqref{E-IF1-ab-gen0fin} reads
\begin{equation}\label{E-IF1-gen}
 \int_M \< q_1\bar A(\beta^{\sharp\top}) +q_2\bar Z,\ \beta^{\sharp}\>\,{\rm d}\vol_a = 0,
\end{equation}
where the constants $q_1$ and $q_2$ are given by
\begin{eqnarray*}
 q_1 \eq -(1+(b^2-\beta(N)^2)\gamma_3)(\hat c\,\rho_1\gamma_1(1+\beta(n)\gamma_1) + \gamma_2(\hat c -\beta(N)\gamma_1)),\\
 q_2 \eq \gamma_1\rho -\hat c\,\rho_1(1+(b^2-\beta(N)^2)\gamma_3)(1+\beta(n)\gamma_1)(\hat c -\beta(N)\gamma_1).
\end{eqnarray*}
\end{thm}

\proof
We calculate
\begin{eqnarray}\label{E-tr-Def}
\nonumber
 \tr\,(\overline{\rm Def}_{\beta^\sharp})_{\,|T\calf}^\top
 \eq\sum\nolimits_{\,i=1}^{\,m}\<\bar\nabla_{b_i}\,\beta^\sharp,\,b_i\>
 =\overline{\Div}\,\beta^\sharp +\beta(\bar Z) - N(\beta(N)),\\
 \<\overline{\rm Def}_{\beta^\sharp}(\beta^{\sharp\top}),\beta^{\sharp\top}\>\eq
 \<\bar\nabla_{\beta^{\sharp\top}}(\beta^{\sharp\top}{+}\beta(N)N),\beta^{\sharp\top}\>
 = \frac12\,\beta^{\sharp\top}(b^2{-}\beta(N)^2)-\beta(N)\<\bar A(\beta^{\sharp\top}),\beta^{\sharp}\>.\qquad
\end{eqnarray}
From Proposition~\ref{L-Dx}, using \eqref{E-tr-Def} and Lemma~\ref{L-zZ}, we get
\[
 \rho\,\|n\|_g\tr A^g =-\tr{\cal A}-\gamma_3\<{\cal A}(\beta^{\sharp\top}),\,\beta^{\sharp}\>,
\]
where
\begin{eqnarray*}
 &&\hskip-9mm \tr{\cal A} = -\rho\,\hat c\tr\bar A -\rho\gamma_1(\,\overline\Div\beta^\sharp +\beta(\bar Z) -N(\beta(N)))
 +m\,n(\rho)/2 +n((b^2-\beta(N)^2)\,\gamma_2)/2 \\
\nonumber
 && -\,\rho\,\beta^{\sharp\top}(\gamma_1)
 +(\rho_0-\rho_1\gamma_1) \big\<\beta(N)\bar\nabla\hat c
 +(\hat c -\beta(N)\gamma_1)(\beta(N)\bar Z -\bar A(\beta^{\sharp\top})) \\
\nonumber
 && -\,(\gamma_1/2)\bar\nabla b^2 -b^2\bar\nabla\gamma_1,\,\beta^{\sharp\top}\big\>
 +\hat c\,\rho_1(1+\beta(n)\gamma_1)\<(\hat c-\beta(N)\gamma_1)\bar Z +\gamma_1\bar A(\beta^{\sharp\top}),\,\beta^{\sharp}\>, \\
 &&\hskip-9mm \<{\cal A}(\beta^{\sharp\top}),\,\beta^{\sharp}\> =
 \rho(\gamma_1\beta(N) -\hat c)\<\bar A(\beta^{\sharp\top}),\,\beta^{\sharp}\>
 -(\rho\,\gamma_1/2)\,\beta^\sharp(b^2-\beta(N)^2) \\
\nonumber
 && +\,(b^2-\beta(N)^2)\,\big[\,n(\rho+(b^2-\beta(N)^2)\,\gamma_2)/2 -\rho\,\beta^{\sharp\top}(\gamma_1)
 +(\rho_0-\rho_1\gamma_1)\big\<\beta(N)\bar\nabla\hat c\\
\nonumber
 && +\,(\hat c -\beta(N)\gamma_1)(\beta(N)\bar Z -\bar A(\beta^{\sharp\top})) -(\gamma_1/2)\bar\nabla b^2
  -b^2\,\bar\nabla\gamma_1,\,\beta^{\sharp\top}\big\>\\
\nonumber
 && +\,\hat c\,\rho_1(1+\beta(n)\gamma_1)
  \<(\hat c-\beta(N)\gamma_1)\bar Z +\gamma_1\bar A(\beta^{\sharp\top}),\beta^{\sharp}\> \big].
\end{eqnarray*}
By \eqref{eq61-g} and \eqref{E-F001vol} for $g$ and $\langle\cdot,\cdot\rangle$,
and reducing terms with factors $b^2-\beta(N)^2$ and $n(\rho)$, we get
\begin{eqnarray*}
\nonumber
 &&\int_M \sigma_g(n)(\rho\,\|n\|_g)^{-1}
 \big\{\,\rho\,\hat c\tr\bar A +\rho\gamma_1(\,\overline\Div\beta^\sharp +\beta(\bar Z) -N(\beta(N))) \\
\nonumber
 && -\,(m-1)\,n(\rho)/2
 +(1+(b^2-\beta(N)^2)\gamma_3)\big[ \rho\,\beta^{\sharp\top}(\gamma_1) -n(\rho+(b^2-\beta(N)^2)\,\gamma_2)/2  \\
\nonumber
 && -\,(\rho_0-\rho_1\gamma_1)\big\<\beta(N)\bar\nabla\hat c
 +\,(\hat c -\beta(N)\gamma_1)(\beta(N)\bar Z -\bar A(\beta^{\sharp\top})) -(\gamma_1/2)\bar\nabla b^2\\
\nonumber
 && -\,b^2\bar\nabla\gamma_1,\,\beta^{\sharp\top}\big\>
 -\hat c\,\rho_1(1+\beta(n)\gamma_1)\<(\hat c-\beta(N)\gamma_1)\bar Z +\gamma_1\bar A(\beta^{\sharp\top}),\,\beta^{\sharp}\>\big]\\
 && +\,\gamma_3\,\rho\big[ (\gamma_1/2)\,\beta^\sharp(b^2-\beta(N)^2)
 +(\hat c-\beta(N)\gamma_1)\<\bar A(\beta^{\sharp\top}),\,\beta^{\sharp}\> \big]\big\}\,{\rm d}\vol_a = 0.
\end{eqnarray*}
The above, equality
\[
 f\,\overline\Div\,\beta^\sharp = \overline\Div\,(f\,\beta^\sharp) -\beta^\sharp(f)
\]
with $f=\sigma_g(n)\,\gamma_1/\|n\|_g$
and the Divergence Theorem yield \eqref{E-IF1-ab-gen0fin}.
For $\beta=0$, we have $\hat c=1$; hence, \eqref{E-IF1-ab-gen0fin} reduces to \eqref{eq61-g} for $a$.
\qed

\begin{example}\rm
(i)~For Randers metric, by Theorem~\ref{T-1-1} we get \cite{rw4}:
\begin{equation*}
 \int_M (c\,\hat c)^{\frac{m}2}c^{-2}\beta(N)\big\{\,\frac12\,N(c^2)
  +\<\bar A(\beta^{\sharp\top}),\beta^{\sharp}\>+c\,\beta(\bar Z)\,\big\}\,{\rm d}\vol_a =0\,,
\end{equation*}
where $c=\hat c-\beta(N)$.
If $b$ and $\beta(N)\ne0$ are constant then $q_1=\hat c\,c(c-\hat c),\,q_2=\hat c(c-\hat c)$ and
\[
 \int_M \<\bar A(\beta^{\sharp\top})+c\bar Z,\,\beta^{\sharp}\>\,{\rm d}\vol_a =0.
\]

(ii)~For Kropina metric, if $b$ and $\beta(N)$ are constant then
\begin{eqnarray*}
 \tr A^g \eq \beta(n)\,\hat c\tr\bar A +\frac{4{\hat c}^{\,2} \beta(n)^2 -2 \beta(n)^2 -\beta(N)^2}{4\,\beta(n)^2}\,\beta(\bar Z)
  +\frac{\beta(N)}{4\,\beta(n)^2}\,\<\bar A(\beta^{\sharp\top}),\beta^{\sharp\top}\>.
\end{eqnarray*}
By Theorem~\ref{T-1-1} we then obtain
\begin{equation*}
 \int_M\<\beta(N)\,\bar A(\beta^{\sharp\top}) +(2\,(2\,{\hat c}^{\,2} -1)\,\beta(n)^2-\beta(N)^2)\,\bar Z,\,\beta^{\sharp\top}\>\,{\rm d}\vol_a =0\,.
\end{equation*}

(iii)~The following application of Theorem~\ref{T-1-1}, when $b$ and $\beta(N)$ are constant,
seems to be inte\-resting.
Let $\bar Z=0$, $q_1\ne0$ and $\alpha$-unit vector field $X\in\Gamma(T\calf)$ be an eigenvector of $\bar{A}$ with an eigenvalue
$\lambda:M\setminus\Sigma\to{\mathbb R}$.
Then $\beta^\sharp=\eps' X+\eps N$, where $\eps =\const\in(0,b_0)$ and $\eps'=\const\in(0,\sqrt{1-\eps^2})$,
obeys \eqref{E-IF1-gen}. Thus, we~get $\int_M \lambda\ {\rm d}\vol_a =0$.
Consequently, either $\lambda\equiv0$ on $M$ or $\lambda(x)\,\lambda(y)<0$ for some points $x$ and $y$ of~$M$.
Furthermore, this implies the Reeb integral formula \eqref{eq61-g} for $a$:
\[
 \int_M (\tr\bar A)\,{\rm d}\vol_a = \sum\nolimits_{\,i}\int_M \lambda_i\,{\rm d}\vol_a = 0.
\]
\end{example}


\begin{thebibliography}{999.}%

\bibitem{r25}
K. Andrzejewski, V. Rovenski, P. Walczak,  {Integral formulas in foliations theory}, 73--82, in
``Geometry \& its Applications", Springer Proc. in Math. \& Statistics, {72}, Springer~(2014).

\bibitem{lw2}
Lu\.{z}y\'{n}czyk M. and Walczak P.
{New integral formulae for two complementary orthogonal distributions on Riemannian manifolds},
Ann. Glob. Anal. Geom. 48 (2015), 195--209.

\bibitem{ma92}
 Matsumoto M. {Theory of Finsler spaces with $(\alpha,\beta)$-metric}, Reports on math. physics, 31(1) (1992), 43--83.

\bibitem{re}
 Reeb G. {Sur la courbure moyenne des vari\'{e}t\'{e}s int\'{e}grales d'une \'{e}quation de Pfaff} $\omega=0$.
 C. R. Acad. Sci. Paris {\bf 231}, 101--102 (1950)

\bibitem{rw1}
 Rovenski V. and Walczak P. \textit{Topics in extrinsic geometry of codimension-one foliations}, Springer, 2011.

\bibitem{rw2}
 Rovenski V. and Walczak P. {Integral formulae on foliated symmetric spaces}, Math. Ann. {352} (2012), 223--237.

\bibitem{rw3}
Rovenski V. and Walczak P. {Integral formulae for codimension-one foliated Finsler spaces},
Balkan J. of Geometry and Its Appl. 21, No. 1 (2016), 76--102 (see ArXiv:1602.00610).

\bibitem{rw4}
 Rovenski V. and Walczak P. Integral formulae for codimension-one foliated Randers spaces,
Publ. Math. Debrecen, vol. 91, 2017.

\bibitem{shsh2016}
 Shen Y.-B. and Shen Z. \textit{Introduction to modern Finsler geometry}, World Scientific, 2016.

\bibitem{sh2}
 Shen Z. \textit{Lectures on Finsler geometry}, World Scientific Publishers, 2001.

\bibitem{yz}
 Yu C. and Zhu H. On a new class of Finsler metrics, Diff. Geom. and its Appl. 29 (2011), 244--254.

\end{thebibliography}
\end{document}